\documentclass[12pt]{article}
\usepackage{amsmath,amssymb,amsthm}

\makeatletter
\@addtoreset{equation}{section}

\title{On $\Z_2$-twisted representation of vertex
operator superalgebras and the Ising model SVOA}

\author{
  Hiroshi Yamauchi
  \vsb\\
  \it{\small Graduate School of Mathematics,}
  \vsb\\
  \it{\small University of Tsukuba, Ibaraki 305-8571, Japan}
  \vsb\\
  {\small e-mail:} \sf{\small  hirocci@math.tsukuba.ac.jp}
}

\date{}

\newcommand{\gs}[1]{\textbf{#1}}

\newcommand{\fr}[2]{\frac{#1}{#2}}
\newcommand{\dfr}[2]{\dfrac{#1}{#2}}
\newcommand{\tfr}[2]{\tfrac{#1}{#2}}
\newcommand{\cd}{\cdot}
\newcommand{\cds}{\cdots}
\newcommand{\dsum}{\displaystyle \sum}

\renewcommand{\l}{\left}
\renewcommand{\r}{\right}

\newcommand{\vsv}{\vspace{5mm}}
\newcommand{\vsb}{\vspace{2mm}}

\newcommand{\q}{\quad}
\newcommand{\qq}{\qquad}

\newcommand{\maru}[1]{{\ooalign{\hfil#1\/\hfil\crcr
\raise.167ex\hbox{\mathhexbox20D}}}}

\newcommand{\ruby}[2]{%
 \leavevmode
 \setbox0=\hbox{#1}%
 \setbox1=\hbox{\tiny #2}%
 \ifdim\wd0>\wd1 \dimen0=\wd0 \else \dimen0=\wd1 \fi
 \hbox{%
   \kanjiskip=0pt plus 2fil
   \xkanjiskip=0pt plus 2fil
   \vbox{%
     \hbox to \dimen0{%
       \tiny \hfil#2\hfil}%
     \nointerlineskip
     \hbox to \dimen0{\mathstrut\hfil#1\hfil}}}}

\newcommand{\lfm}{\langle} 
\newcommand{\rfm}{\rangle} 

\newcommand{\abs}[1]{\lvert{#1}\rvert}

\DeclareMathOperator*{\tensor}{\otimes}

\newcommand{\Z}{\mathbb{Z}}
\newcommand{\C}{\mathbb{C}}
\newcommand{\R}{\mathbb{R}}
\newcommand{\N}{\mathbb{N}}

\newcommand{\res}{\mathrm{Res}}
\newcommand{\End}{\mathrm{End}}

\newcommand{\aut}{\mathrm{Aut}}
\newcommand{\wt}{\mathrm{wt}}

\renewcommand{\hom}{\mathrm{Hom}}

\newcommand{\id}{\mathrm{id}}

\newcommand{\pii}{\pi i}

\newcommand{\w}{\omega}
\newcommand{\vacuum}{\mathrm{1\hspace{-3.2pt}l}}
\newcommand{\vac}{\vacuum}

\newcommand{\ising}{L(\tfr{1}{2}, 0)}
\newcommand{\isingodd}{L(\tfr{1}{2},\tfr{1}{2})}
\newcommand{\isingmod}{L(\tfr{1}{2},\tfr{1}{16})}
\newcommand{\g}{\mathfrak{g}}
\newcommand{\simmapsto}{\stackrel{\sim}{\mapsto}}
\newcommand{\hf}{\frac{1}{2}}
\newcommand{\thf}{\tfrac{1}{2}}
\newcommand{\pf}{\gs{Proof:}\q}
\newcommand{\VB}{V\!\! B^\natural}
\newcommand{\B}{\mathbb{B}}
\newcommand{\M}{\mathbb{M}}
\newcommand{\nat}{\natural}
\newcommand{\VBT}{\VB_{\mathrm{tw}}}
\newcommand{\irr}{\mathrm{Irr}}

\theoremstyle{plain}
\newtheorem{thm}{Theorem}[section]
\newtheorem{prop}[thm]{Proposition}
\newtheorem{lem}[thm]{Lemma}

\theoremstyle{definition}
\newtheorem{df}[thm]{Definition}

\theoremstyle{remark}
\newtheorem{rem}[thm]{Remark}

\begin{document}

\baselineskip 6mm

\maketitle

\begin{abstract}
  We investigate a general theory of the $\Z_2$-twisted
  representations of vertex operator superalgebras. 
  Certain one-to-one correspondence theorems are  established.
  We also give an explicit realization of the Ising model SVOA and its
  $\Z_2$-twisted modules.
  As an application, we obtain the Gerald H\"{o}hn's Babymonster SVOA 
  $\VB$ and its $\Z_2$-twisted module $\VBT$ from the moonshine
  VOA $V^\nat$ by cutting off the Ising models. It is also shown in
  this paper that $\aut (\VB )$ is finite.
\end{abstract}

\section{Introduction} 

In the theory of vertex operator algebras (VOAs), we sometimes notice
that it makes the theory simpler to use some representations of
vertex operator superalgebras (SVOAs for short) instead of VOAs.
For example, as one can see in \cite{M1}-\cite{M3}, in the
representation theory of the Ising model VOA $\ising$ it seems more
natural that we should treat the theory in the view point of the SVOA
$\ising\oplus\isingodd$, which is the most interested object in this
paper.
For the theory of VOAs, we have many remarkable results on the
fundamental representation theory, so-called  Zhu's theory
(cf. \cite{Z} \cite{FZ}  \cite{Li2} \cite{DLM1} \cite{DLM2} \cite{Y})
and some of them are extended to that for SVOAs (cf. \cite{KW}). 
The {\it Zhu algebra} $A(V)$ is an associative algebra associated to
every VOA $V$. It is well-known that there exists a one-to-one
correspondence between the category of irreducible $V$-modules and
that of irreducible $A(V)$-modules (cf. \cite{Z}, \cite{DLM1}).
Since every SVOA has a canonical involution, we can think of the
$\Z_2$-twisted representations of SVOAs and the above theory should be 
naturally extended to $\Z_2$-twisted representations for SVOAs.

To start the investigation of SVOAs, the Ising model SVOA
$\ising\oplus\isingodd$ will be a good example.  
It is one of the smallest SVOAs and a fundamental object in the theory
of rational VOAs. 
In spite of its simplicity, there are many applications and deep 
theories. 
The fusion algebra of the Ising model VOA $\ising$ has been determined
and many related results are obtained (cf. \cite{DMZ}, \cite{DGH},
\cite{M1}-\cite{M3}, etc). 
In particular, in \cite{M3} the moonshine vertex operator algebra
$V^\natural$ is reconstructed by Miyamoto by using the Ising models. 
In the Miyamoto's theory, some central extensions of $2$-groups play
an important role. 
It seems that the appearance of $2$-groups reflects the
$\Z_2$-gradation of the structure of SVOAs. 
Since there are two ways to represent the $\Z_2$-gradation of SVOA,
namely, we can consider the $\Z_2$-graded representation and
the $\Z_2$-twisted representation, the latter shall be as important as 
the former. 
One can also find the importance of the $\Z_2$-twisted representation
from another point of view.
In Frenkel-Lepowsky-Muerman's $\Z_2$-orbifold construction \cite{FLM},
a $\Z_2$-twisted representation of the Leech lattice VOA is used to
construct the moonshine VOA.  
Through the two different constructions of the moonshine VOA, we come
to believe that the $\Z_2$-twisted representation of SVOA will be an
essential concept in the theory of SVOAs.
In this respect, we introduce the $\Z_2$-twisted Zhu algebras
and investigate the Ising model SVOA.
We will give an explicit realization of the Ising model SVOA and its
all $\Z_2$-twisted representations. 
Namely, we will give a realization of the all unitary representation
of the Virasoro algebra of central charge $\hf$ in terms of local
and $\Z_2$-twisted local systems. 
It enable us to compute all type of intertwining operators in
fermionic literature, which is simpler than the known bosonic
literature.
Based on our realization, we will determine the structure of
$\Z_2$-twisted Zhu algebra for the Ising model SVOA.
As an application, we propose a method to construct some
SVOAs and its $\Z_2$-twisted representations using the Ising models.
Applying our method on the moonshine VOA $V^\nat$, we obtain the
Gerald H\"{o}hn's Babymonster SVOA $\VB$ and its irreducible
$\Z_2$-twisted representation $\VBT$.

We organize this paper in the following way.
In Sec. 2 we generalize Kac-Wang's Zhu algebra associated to SVOA to
the twisted case and then we establish a bijective correspondence which
is well-known in the theory of VOAs. In the proof of the associativity 
of the Zhu algebra, by making use of the commutativity and
associativity of an SVA, we can make the proof much simpler.
We also present Frenkel-Zhu's bimodules for
$\Z_2$-twisted Zhu algebras.
The fusion rules for $\Z_2$-twisted modules are described in terms of
our bimodules for $\Z_2$-twisted Zhu algebras.

In Sec. 3 we present some basic facts on SVOA.
Invariant bilinear forms for SVOA and $\Z_2$-conjugacy of the modules
are treated.

In Sec. 4 we present a construction of the Ising models.
Using the ideas of the local system and twisted local system in
\cite{Li1}, we attain an SVOA structure in $M=\ising\oplus\isingodd$
and its $\Z_2$-twisted representation in $\isingmod$ in Subsec 4.2.  
Although some ingredients in Subsec 4.1-4.2 are already known in some
papers (cf. \cite{FRW},\cite{FRW2},\cite{La}, etc), we will pick up
and repeat the necessary parts with the suitable modifications because
fermionic construction has one merit such that it enable us to compute
intertwining operators of all type for the Ising models including the
twisted part.
In Subsec. 4.3 we determine $\Z_2$-twisted Zhu algebras for the Ising
model SVOA and classify all irreducible $\Z_2$-twisted
representations based on our realization given previously. 
At the last of this section, we consider an application of the Ising
models. 
We present a method to construct an SVOA from a VOA containing the
Ising models. 
This method is already shown in \cite{H}.
However, it is worth to present our method since it works in simpler
situation. 

In Sec. 5 we present a construction of the Gerald H\"{o}hn's
Babymonster SVOA and its $\Z_2$-twisted module from the moonshine VOA.
It is a simple SVOA whose automorphism group contains Fischer's Baby
monster sporadic simple group.
Using some methods from the Quantum Galois theory, we will show that
the full automorphism group of the Babymonster SVOA is finite.

\section{$\Z_2$-twisted Zhu theory for SVOAs} 

\subsection{$\Z_2$-twisted representations} 

Let $(V,Y,\vac,\w )$ be an SVOA, where $V$ has a $\Z_2$-grading $V=V^0 
\oplus V^1$. We assume that $V$ has a half-integer grading: $V^0= 
\oplus_{n\in\Z}V^0_n$ and $V^1=\oplus_{n\in\Z}V^1_{n+1/2}$, where
$V^i_s=\{ v\in V | L_0 v=s v\}$, $i\in \Z_2$, $s\in \fr{1}{2}\Z$.
We also assume that each $L_0$-weight space $V_s$ is of finite
dimension.
Let $p$ and $q$ be the parity functions defined by $p(a,b)=1$ for
$a,b\in V^1$ and $p(a,b)=0$ for $\Z_2$-homogeneous $a,b\in V$ of the
other cases, and $q(a)=i$ for $a\in V^i$, respectively.
In this paper, we will treat two distinct representations of
$V$. Even though these definitions are already shown in
\cite{Li1}, we repeat them for convenience.

\begin{df}
  A $\Z_2$-graded $V$-module is a pair $(M, Y_M)$ consisting of a
  $\Z_2$-graded vector space $M=M^0\oplus M^1$ which comes from a
  $\fr{1}{2}\N$-grading $M^i=\oplus_{n\in\N +\fr{i}{2}}^\infty M^i(n)$
  $(i=0,1)$ and a linear map $Y_M(\cd\, ,z)$ from $V$ to $\End (M)
  [[z,z^{-1}]]$ satisfying the following conditions: 
  \vsb\\
  \q $1^\circ$\ For any $a\in V$, $v\in M$, $a_nv=0$ for $n$
     sufficiently large;
  \vsb\\
  \q $2^\circ$\ $Y_M(\vac ,z)=\id_M$;
  \vsb\\
  \q $3^\circ$\ $a_n M(s)\subset M(s+\wt (a) -n-1)$;
  \vsb\\
  \q $4^\circ$\ For any $\Z_2$-homogeneous $a,b\in V$, the following
      Jacobi identity holds:

  \begin{equation}
    \begin{array}{c}
    z_0^{-1}\delta\l(\dfr{z_1-z_2}{z_0}\r) Y_M(a,z_1)Y_M(b,z_2)
    -(-1)^{p(a,b)}z_0^{-1}\delta\l(\dfr{-z_2+z_1}{z_0}\r)
      Y_M(b,z_2)Y_M(a,z_1)
    \vsb\\
    =z_2^{-1}\delta\l(\dfr{z_1-z_0}{z_2}\r) Y_M(Y_V(a,z_0)b,z_2) .
    \end{array}
  \end{equation}
\end{df}

\begin{df}
  A $\Z_2$-twisted $V$-module is a pair $(M,Y_M)$ consisting of an
  $\N$-graded vector space $M=\oplus_{n\in\N} M(n)$ and a linear map
  $Y_M(\cd\, ,z)$ from $V$ to $\End (M) 
  [[z^{\fr{1}{2}}, z^{-\fr{1}{2}}]]$ satisfying the following
  conditions:  
  \vsb\\
  \begin{tabular}{ll}
    $1^\circ$ & For $a\in V^i$, the module vertex operator has the
      shape $Y_M(a,z) = \dsum_{n\in\fr{i}{2}+\Z} 
      a_{n}z^{-n-1}$;
    \vsb\\
    $2^\circ$ & For any $a\in V$, $v\in M$, $a_nv=0$ for $n\in
      \fr{1}{2}\Z$  sufficiently large;
    \vsb\\
    $3^\circ$ & $Y_M(\vac ,z)=\id_M$;
    \vsb\\
    $4^\circ$ & $a_n M(s)\subset M(s+\wt (a) -n-1)$;
    \vsb\\
    $5^\circ$ & For any $\Z_2$-homogeneous $a,b\in V$, the following
      $\Z_2$-twisted Jacobi identity holds:
  \end{tabular}
  \begin{equation}
    \begin{array}{c}
    z_0^{-1}\delta\l(\dfr{z_1-z_2}{z_0}\r) Y_M(a,z_1)Y_M(b,z_2)
    -(-1)^{p(a,b)}z_0^{-1}\delta\l(\dfr{-z_2+z_1}{z_0}\r)
      Y_M(b,z_2)Y_M(a,z_1)
    \vsb\\
    = z_2^{-1}\delta\l(\dfr{z_1-z_0}{z_2}\r) 
      \l(\dfr{z_1-z_0}{z_2}\r)^{-\fr{q(a)}{2}} Y_M(Y_V(a,z_0)b,z_2) .
    \end{array}
  \end{equation}
\end{df}

\subsection{$\Z_2$-twisted Zhu algebras} 

As in the case of VOAs, we can define the Zhu algebras for SVOAs.
Since the representations of Zhu algebras correspond to those of
original SVOAs, we can introduce two type of Zhu algebras for
corresponding type of the representations. 
First, we consider Zhu algebras for the $\Z_2$-graded representations. 
The following definition is due to Kac-Wang \cite{KW}.

\begin{df}
  We define the bilinear maps $*: V\tensor V\to V$, $\circ : V\tensor
  V\to V$ as follows.
  $$
  \begin{array}{lll}
     a*b & := & 
       \begin{cases}
         \ \res_z Y(a,z)\dfr{(1+z)^{\wt (a)}}{z}b 
         & \text{if}\ a\in V^0,
         \vsb\\
         \ 0 & \text{if}\ a\in V^1,
       \end{cases}
     \vsb\\
     a\circ b & := & 
       \begin{cases}
         \ \res_z Y(a,z)\dfr{(1+z)^{\wt (a)}}{z^2} b
         & \text{for}\ a\in V^0,
         \vsb\\
         \ \res_z Y(a,z)\dfr{(1+z)^{\wt (a)-\fr{1}{2}}}{z}b
         & \text{for}\ a\in V^1.
       \end{cases}
  \end{array}
  $$
  Extend to $V\tensor V$ linearly, denote by $O(V)\subset V$ the
  linear span of elements of the form $a\circ b$, and by $A(V)$ the
  quotient space $V/O(V)$.
\end{df}

\begin{rem}
  It follows from the definition that $a\circ \vac =a$ for $a\in V^1$ 
  so that $V^1\subset O(V)$. Notice that $O(V^0)\subset
  O(V)$ where $O(V^0)$ is the kernel of the Zhu algebra $A(V^0)$ for 
  a VOA $V^0$. Therefore, $A(V)$ is a quotient algebra of $A(V^0)$.
\end{rem}

In \cite{KW} one can find the followings.

\begin{thm} (Theorem 1.1, 1.2, 1.3 in \cite{KW})
  \vsb\\
  (1)\ $O(V)$ is a two-sided ideal of $V$ under the
    multiplication $*$. Moreover, the quotient algebra $(A(V), *)$ is 
    associative.
  \vsb\\
  (2)\ $\vac +O(V)$ is the unit element of $A(V)$ and $\w +O(V)$ is
    in the center of $A(V)$.
  \vsb\\
  (3)\ Let $M=\oplus_{n\in\fr{1}{2}\N}M_n$ be a $\Z_2$-graded
    $V$-module. 
    Then the top level $M_0$ is an $A(V)$-module via $a+O(V)\mapsto
    o(a)=a_{\wt (a)-1}$.
  \vsb\\
  (4)\ Given an $A(V)$-module $(W,\pi )$, there exists a
    $\Z_2$-graded $V$-module $M=\oplus_{n\in \fr{1}{2}\N} M_n$ such
    that the $A(V)$-modules $M_0$ and $W$ are isomorphic. 
    Moreover, this gives a bijective correspondence between the set of
    irreducible $A(V)$-modules and the set of irreducible
    $\Z_2$-graded $V$-modules. 
\end{thm}

We can define another Zhu algebra, the one for $\Z_2$-twisted
representations as follow.

\begin{df} $(\Z_2$-twisted Zhu algebra$)$
  \vsb\\
  For $L_0$-homogeneous $a,b\in V$, define the bilinear maps $*_t:
  V\tensor V\to V$, $\circ_t : V\tensor V\to V$ as follows.
  $$
  \begin{array}{lll}
     a*_t b & := & \res_z Y(a,z)\dfr{(1+z)^{\wt (a)}}{z}b ,
     \vsb\\
     a\circ_t b & := & \res_z Y(a,z)\dfr{(1+z)^{\wt (a)}}{z^2} b.
  \end{array}
  $$
  Extend to $V\tensor V$ linearly, denote by $O_t(V)\subset V$ the
  linear span of elements of the form $a\circ_t b$, and by $A_t(V)$
  the quotient space $V/O_t(V)$.
\end{df}

We will show

\begin{thm}\label{1:1}
  (1)\ The quotient space $(A_t(V),*_t)$ is a $\Z_2$-graded
  associative algebra with unit element $\vac +O_t(V)$.
  \vsb\\
  (2)\ Let $M=\oplus_{n\in\N} M(n)$ be a $\Z_2$-twisted
  $V$-modules. Then the top level $M(0)$ of $M$ is an $A_t(V)$-module
  under the action $a+O_t(V)\mapsto o(a)=a_{\wt (a)-1}$ for
  homogeneous $a\in V$.
  \vsb\\
  (3)\ Let $(W,\pi )$ be an irreducible $A_t(V)$-module.
  Then there exists an irreducible $\Z_2$-twisted $V$-module
  $M=\oplus_{n\in\N} M(n)$ such that the top level $M(0)$ of $M$ 
  is isomorphic to $W$ as $A_t(V)$-modules. 
\end{thm}  

\pf
Since we can prove the above statements by the same argument used in
\cite{DLM1}, we give a slight one here.
To prove (1), using the associativity of a vertex algebra, we can make 
proof very simpler.
Let $A$ be an associative algebra $\C [t^\fr{1}{2},t^{-\fr{1}{2}}]$.
Then $(A,Y_A,1,\fr{d}{dt})$ is a vertex algebra on which the vertex
operator is defined by $Y(t^s,z)=e^{z\fr{d}{dt}}t^s=(t+z)^s$.
$A$ has a vertex algebra grading $A = \oplus_{s\in\fr{1}{2}\Z} A_s$
defined by $\wt (t^s):= -s$. 
Let $\hat{V}:=A \tensor_\C V$ be a tensor product of vertex
(super)algebras. It is clear that $(\hat{V}, Y_A\tensor Y_V, 1\tensor
\vac, \fr{d}{dt}\tensor 1+1\tensor L_{-1})$ is a vertex superalgebra
and $\hat{V}$ also carries a $\fr{1}{2}\Z$-grading $\hat{V} = 
\oplus_{n \in \fr{1}{2}\Z} \hat{V}(n)$, where
$\hat{V}(n)=\oplus_{i+j=n} A_i\tensor V_j$. 
Define a linear isomorphism $\theta$ of $\hat{V}$ by
$\theta (t^s\tensor a) := (-1)^{q(a)} e^{-2\pii s} \cd t^s \tensor a$.  
Then one can easily check that $\theta$ is an isomorphism of vertex
superalgebra $\hat{V}$. Take $\theta$-invariants of $\hat{V}$ and
denote it by $\hat{V}^\theta$. It is obvious that $\hat{V}^\theta$
admits a $\Z$-grading decomposition
$$
\begin{array}{ll}
  \hat{V}^\theta = \bigoplus_{n\in\Z} \hat{V}^\theta(n),
  & 
  \hat{V}^\theta (n) = \bigoplus_{i+j=n} A_i\tensor V_j.
\end{array}
$$
A subspace $\hat{V}^\theta (0)$ is an algebra with unit element
$1\tensor \vac$ under the multiplication $X \cd Y:=X \bullet_{-1}Y$
for $X,Y\in \hat{V}^\theta(0)$, where $\bullet_n$ denotes the $n$-th
product in $\hat{V}$.
As a linear space, $\hat{V}^\theta (0)=\oplus_{s\in \fr{1}{2}\Z_+}
A_{-s}\tensor V_s$ is isomorphic to $V$ under the mapping 
$t^{\wt (a)}\tensor a\mapsto a$, so we can identify them. 
Then $O_t(V)$ is isomorphic to $\hat{V}^\theta (0)\bullet_{-2}
\hat{V}^\theta(-1)$, 
i.e.
$$
\begin{array}{lll}
  a\circ_t b &=& \res_z Y(a,z)b\dfr{(1+z)^{\wt (a)}}{z^2}
  \vsb\\
  &\simmapsto& \res_z Y(a,z)\dfr{(t+z)^{\wt (a)}}{z^2}\cd 
  (t+z)^{\wt (b) +1}b 
  \vsb\\
  &=& (t^{\wt (a)}\tensor a) \bullet_{-2}(t^{\wt (b) +1}\tensor b) ,
\end{array}
$$
under the identification. 
Therefore, $A_t(V)$ is linearly isomorphic to
$\hat{V}^\theta(0)/K$, $K = \hat{V}^\theta(0) \bullet_{-2}
\hat{V}^\theta (-1)$.
Furthermore, under the identification, the product $a*_t b$ in
$A_t(V)$ corresponds to 
$(t^{\wt (a)}\tensor a) \bullet_{-1} (t^{\wt (b)+1}\tensor b)$.  
So to prove (1), we should show that $(\hat{V}^\theta (0)/K,
\bullet_{-1})$ is an associative algebra. 
To show this, we need the following simple lemma which is a direct
consequence of the definition of the tensor product of SVAs. 

\begin{lem}\label{affinezhu}
  (Lemma 2.1.2 in \cite{Z})\ 
  For any $n\in\Z$, $i,j\in\N$, we have the followings.
  $$
  \begin{array}{c}
    \hat{V}^\theta(-i) \bullet_{-2-i} \hat{V}^\theta(n)
      \subset \hat{V}^\theta(0) \bullet_{-2} \hat{V}^\theta(n),
    \vsb\\
    \hat{V}^\theta(0) \bullet_{-2-j} \hat{V}^\theta(n-j)
      \subset \hat{V}^\theta(0) \bullet_{-2}\hat{V}^\theta(n).
  \end{array}
  $$
  In particular, $\hat{V}^\theta(-i) \bullet_{-2-i-j} \hat{V}^\theta
  (-1-j) \subset K$ holds for any $i,j\in\N$. 
\end{lem}

Let $E,F,G\in \hat{V}^\theta(0)$ and $H\in \hat{V}^\theta(-1)$.
To prove that $K$ is an two-sided ideal of $\hat{V}^\theta(0)$, 
it is suffice to show that $E_{-1}F_{-2}H$ and $(F_{-2}H)_{-1}E$ 
belong to $K$. Note that $\hat{V}^\theta$ is a sub SVA of $\hat{V}$.
The commutativity and associativity of SVA and Lemma
\ref{affinezhu} lead
$$
\begin{array}{ll}
  E_{-1}F_{-2}H 
  &= (-1)^{p(E,F)} F_{-2}E_{-1}H +[E_{-1},F_{-2}]H
  \vsb\\
  &= (-1)^{p(E,F)} F_{-2}E_{-1}H +\dsum_{i=0}^\infty 
    \dbinom{-1}{i} (E_iF)_{-3-i}H
    \in K,
  \vsb\\
  (F_{-2}H)_{-1}E 
  &= \dsum_{i=0}^\infty (-1)^i\dbinom{-2}{i}\l\{ 
    F_{-2-i}H_{-1+i} -(-1)^{p(F,H)}H_{-3-i}F_i\r\} E \in K.
\end{array}
$$
Thus $K$ is a two-sided ideal of $\hat{V}^\theta (0)$.
Similarly, we have
$$
\begin{array}{l}
  E_{-1}F_{-1}G -(E_{-1}F)_{-1}G
  \vsb\\
  = E_{-1}F_{-1}G -\dsum_{i=0}^\infty (-1)^i \dbinom{-1}{i}
    \l\{ E_{-1-i}F_{-1+i}-(-1)^{p(E,F)-1}F_{-2-i}E_i\r\} G
  \vsb\\
  = - \dsum_{i=1}^\infty E_{-1-i}F_{-1+i}G 
     + \dsum_{i=0}^\infty (-1)^{p(E,F)-1}F_{-2-i}E_iG
     \in W.
\end{array}
$$
Therefore $\hat{V}^\theta (0)/K$ is associative.
Since $O_t(V)$ has a natural $\Z_2$-grading induced from $V$,
$A_t(V)=V/O_t(V)$ also has a natural $\Z_2$-grading induced from
$V$. This proves (1). 

(2) is similar to Theorem 2.1.2 of \cite{Z} so that we consider (3).
Denote the derivation operator $\fr{d}{dt}\tensor 1+1\tensor L_{-1}$
on $\hat{V}$ by $D$. Then $(\hat{V}^\theta /D\hat{V}^\theta,
\bullet_0)$ becomes a Lie superalgebra (cf. \cite{Bo}).
Denote this Lie superalgebra by $\g (V)$. Defining a $\Z$-gradation
on $\g (V)$ by $\deg (t^n\tensor a+D\hat{V}^\theta ):= \wt (a)-n-1$, we
obtain a $\Z$-grading decomposition $\g (V) = \oplus_{n\in\Z} \g (V)_n$
of Lie superalgebra and using the grading we have a triangular
decomposition of a Lie superalgebra $\g (V)$ as follow. 
$$
  \g (V)=\g (V)^-\oplus \g (V)^0\oplus \g (V)^+,
$$
where $\g (V)^\pm =\oplus_{\pm n> 0}\g (V)^n$ and $\g (V)^0=\g
(V)_0$.
Let $U(\g (V))$ be the universal enveloping algebra for $\g (V)$.
Replacing $U(V[g])$ in \cite{DLM1} by $U(\g (V))$, we can prove (3)
by exactly the same argument as that in \cite{DLM1}.
\qed

\begin{rem}
  One can show that the image of $\w$ in $A_t(V)$ is in the center of
  $A_t(V)$ by a direct calculation. But our ``affinization argument''
  can't cover this proof.
\end{rem}

\begin{rem}\label{tau}
  Since $V$ has a canonical involution $\sigma$ which is identical on 
  $V^0$ and acts as $-1$ on $V^1$, we can consider a $\sigma$-twisted
  $V$-module, which is exactly a $\Z_2$-twisted $V$-module in our
  notation. 
  Contrary to the terminology 'twisted', it seems that 
  $\Z_2$-twisted Zhu algebras for SVOAs correspond to non-twisted
  Zhu algebras for VOAs and non-twisted ones for SVOAs
  correspond to twisted ones for VOAs.
\end{rem}

\subsection{Frenkel-Zhu's bimodules and intertwining operators} 

In this subsection we define a bimodule $A_t(U)$ of $A_t(V)$ for every
$\Z_2$-graded $V$-module $U$ as a generalization of \cite{FZ} and
\cite{KW}.
Then we give a description of the fusion rules among $\Z_2$-twisted
modules in terms of $A_t(U)$.

Let $M^1$, $M^2$, $M^3$ be irreducible $\Z_2$-graded $V$-modules with 
$\hf \N$-grading $M^i=\oplus_{n\in\hf \N}$ $M^i(n)$ and
$\Z_2$-grading $(M^i)^r=\oplus_{n\in \N}M^i(n+\fr{r}{2})$ for
$i=1,2,3$ and $r=0,1$, respectively.
It is known that all irreducible $V$-modules admit the $L_0$-weight
space decomposition so that  we can find some $h_i\in \C$ such that 
$M^i(n)=M^i_{n+h_i}$ for $i=1,2,3$, respectively, where $X_s$ denotes
the $L_0$-weight space of $X$ with weight $s\in\C$.

\begin{df}
  Under the above setting, a $\Z_2$-graded intertwining operator of
  type $\binom{M^3}{M^1\ M^2}$ is a linear map
  $$
    I(\cd ,z) : u\in M^1 \mapsto I(u,z)=\dsum_{s\in\C} u_s
    z^{-s-1} \in \hom_\C (M^2,M^3)\{ z\}
  $$
  satisfying the following conditions:
  \vsb\\
  $1^\circ$\q For every $v\in M^1$, 
    $I(u,z)\in \hom_\C (M^2, M^3)[[z,z^{-1}]] z^{-h_1-h_2+h_3}$;
  \vsb\\
  $2^\circ$\q For any $s\in \C$, $u\in M^1$ and  $v\in M^2$,
    $u_{s+N}v=0$ for sufficiently large $N\in \N$;
  \vsb\\
  $3^\circ$\q $L_{-1}$-derivation: $I(L_{-1}u,z)=\dfr{d}{dz} I(u,z)$; 
  \vsb\\
  $4^\circ$\q For $\Z_2$-homogeneous $a\in V$ and $u\in M^1$, the
  following Jacobi identity holds:
  $$
  \begin{array}{c}
    z_0^{-1} \delta\l(\dfr{z_1-z_2}{z_0}\r) Y_{M^3}(a,z_1) I(u,z_2)
      -(-1)^{p(a,u)} z_0^{-1} \delta\l(\dfr{-z_2+z_1}{z_0}\r) I(u,z_2)
      Y_{M^2}(a,z_2)
    \vsb\\
    = z_2^{-1} \delta\l(\dfr{z_1-z_0}{z_2}\r) I(Y_{M^1}(a,z_0)u,z_2) ,
  \end{array}
  $$
  where the pairing $p(\cd ,\cd ): V\times U\to \Z_2$ is understood
  appropriately. 
\end{df}

Frenkel-Zhu's bimodules for SVOAs are introduced by Kac-Wang in
\cite{KW}. 

\begin{df}
  For a $\Z_2$-graded $V$-module $U$, we define bilinear operations
  $a\circ u$, $a*u$ and $u*a$, for $a\in V$ homogeneous and $u\in U$,
  as follows 
  $$
  \begin{array}{l}
    a\circ u:= \res_z Y(a,z)\dfr{(1+z)^{\wt (a)-\fr{r}{2}}}{z^{2-r}}
    u,\q \text{for}\ a\in V^r,
    \vsb\\
    a*u:= \res_z Y(a,z)\dfr{(1+z)^{\wt (a)}}{z} u,\q \text{for}\ a\in
    V^0,
    \vsb\\
    u*a:= \res_z Y(a,z)\dfr{(1+z)^{\wt (a)-1}}{z} u,\q \text{for}\
    a\in V^0,
    \vsb\\
    a*u=u*a=0,\q \text{for}\ a\in V^1
  \end{array}
  $$
  and extend linearly. We also define $O(U)\subset U$ to be the 
  linear span of elements of the form $a\circ u$ and $A(U)$ to be the
  quotient space $U^0/\l( O(U)\cap U^0\r)$.
\end{df}

\begin{rem}
  Our definition of $A(U)$ differs from Kac-Wang's original one.
  Namely, we define $A(U)$ to be a quotient space of $U^0$.
  See \cite{Y} for the validity of this change.
\end{rem}

By definition, we know that $z^{h_1+h_2-h_3} I(\cd, z)\in \hom_\C
(M^2,M^3) [[z,z^{-1}]]$. 
It is convenient to set $I(u,z)=\sum_{n\in \Z}
u_{(n)}z^{-n-1-h_1-h_2+h_3}$ and $\deg (u) := \wt (u)-h_1$ for $u\in
M^1$ to state the following theorems which are shown in \cite{KW}.

\begin{thm} \label{thm.1.4} (Theorem 1.4 in \cite{KW}) 
  $A(U)$ is an $A(V)$-bimodule under the action $*$.
\end{thm}

\begin{thm} \label{thm.1.5} 
  (1)\ (Theorem 1.5 in \cite{KW}) 
  For a $\Z_2$-graded intertwining operator $I(\cd,z)$ of
  type $\binom{M^2}{U\ M^1}$, its zero-mode action $o^I(u):=
  u_{(\deg (u)-1)}$ gives a linear injection from $I\binom{M^2}{U\
  M^1}$ to $\hom_{A(V)}$ $\big( A(U)\tensor_{A(V)}$ $M^1(0)$,
  $M^2(0)\big)$. 
  \vsb\\
  (2)\ (Theorem 2.11 in \cite{Li2})
  Suppose that every $\Z_2$-graded $V$-module is completely
  reducible. Then the linear map $I(\cd ,z) \mapsto o^I$ given in (1)
  defines a linear isomorphism of vector spaces $I\binom{M^2}{U\ M^1}$
  and $\hom_{A(V)}$ $\big( A(U)\tensor_{A(V)}M^1(0)$, $M^2(0)\big)$.
\end{thm}

\begin{rem}
  As pointed out in \cite{Li2}, the assumption on completely
  reducibility in the above statement (2) is necessary.
  For details, see \cite{Li2}. 
\end{rem}

We extend the above results to $\Z_2$-twisted case.
We begin with introducing the notion of $\Z_2$-twisted intertwining
operators.
Let $U$ be an irreducible $\Z_2$-graded $V$-module and 
$W^1$, $W^2$ be irreducible $\Z_2$-twisted $V$-modules with the
weight space decompositions $U=\oplus_{n\in\hf\N} U_{n+h_0}$ and  
$W^i=\oplus_{n\in\N} W^i_{n+h_i}$ for $i=1,2$, respectively.

\begin{df}
  Under the above setting, a $\Z_2$-twisted intertwining operator
  $I(\cd ,z)$ of type $\binom{W^2}{U\ W^1}$ is a linear map
  $$
    I(\cd,z) : u\in U\mapsto I(u,z)=\dsum_{s\in\C} u_s z^{-s-1}
    \in \hom_\C (W^1,W^2)\{ z\} 
  $$
  satisfying the following conditions:
  \vsb\\
  $1^\circ$\q For $u\in U^r$, $I(u,z)\in \hom_\C (W^1, W^2)
    [[z,z^{-1}]] z^{-\fr{r}{2}-h_0-h_1+h_2}$, where $r=0,1$;
  \vsb\\
  $2^\circ$\q For any $s\in \C$, $u\in U$ and $w\in W^1$, $u_{s+N}w=0$ 
    for sufficiently large $N\in \N$;
  \vsb\\
  $3^\circ$\q $L_{-1}$-derivation: $I(L_{-1}u,z)=\dfr{d}{dz}I(u,z)$;
  \vsb\\
  $4^\circ$\q For $\Z_2$-homogeneous $a\in V$ and $u\in U$, the
  following $\Z_2$-twisted Jacobi identity holds:
  $$
  \begin{array}{c}
    z_0^{-1}\delta\l(\dfr{z_1-z_2}{z_0}\r) Y_{W^2}(a,z_1) I(u,z_2)
      -(-1)^{p(a,u)} z_0^{-1} \delta\l(\dfr{-z_2+z_1}{z_0}\r) 
      I(u,z_2) Y_{W^1}(a,z_1)  
    \vsb\\
    = z_2^{-1}\delta\l(\dfr{z_1-z_0}{z_2}\r) \l(\dfr{z_1-z_0}{z_2}
    \r)^{-\fr{q(a)}{2}} I(Y_U(a,z_0)u,z_2) .
  \end{array}
  $$
\end{df}

As we did previously, we set $I(u,z)=\sum_{n\in\Z +\fr{r}{2}} u_{(n)} 
z^{-n-1-h_0-h_1+h_2}$ and $\deg (u):=\wt (u)-h_0$ for $u\in U^r$. 

\begin{rem}
  For a $\Z_2$-twisted $V$-module $W$, the module vertex operator
  $Y_W(\cd ,z)$ is a $\Z_2$-twisted intertwining operator of type 
  $\binom{W}{V\ W}$ by definition.
\end{rem}

We define the following products.

\begin{df}
  For a $\Z_2$-graded $V$-module $U$, we define bilinear operations
  $a\circ_t u$, $a*_t u$ and $u*_t a$, for homogeneous $a\in V$ and
  $u\in U$, as follows
  $$
  \begin{array}{l}
    a\circ_t u:= \res_z Y(a,z) \dfr{(1+z)^{\wt (a)}}{z^2} u,
    \vsb\\
    a*_t u:= \res_z Y(a,z) \dfr{(1+z)^{\wt (a)}}{z} u,
    \vsb\\
    a*_t u:= (-1)^{p(a,u)} \res_z Y(a,z) \dfr{(1+z)^{\wt (a)-1}}{z}u ,
  \end{array}
  $$
  and extend linearly. 
  We also define $O_t(U)\subset U$ to be the linear span of elements
  of the form $a\circ_t u$ and $A_t(U)$ to be the quotient space
  $U/O_t(U)$. 
\end{df}

As an analogy of Theorem \ref{thm.1.4} and \ref{thm.1.5}, we have the 
followings.

\begin{thm}
  (1)\ $A_t(U)$ is an $A_t(V)$-bimodule under the action $*_t$.
  \vsb\\
  (2)\ For a $\Z_2$-twisted intertwining operator $I(\cd ,z)$ of type
    $\binom{W^2}{U\ W^1}$, its zero-mode action $o^I(u):= u_{(\deg
    (u)-1)}$ gives a linear injection from $I\binom{W^2}{U\ W^1}$ to 
    $\hom_{A_t(V)}$ $\big( A_t(U)\tensor_{A_t(V)}$ $W^1(0)$,
    $W^2(0)\big)$. 
  \vsb\\
  (3)\ Suppose that every $\Z_2$-twisted $V$-module is completely
  reducible. Then the linear map $I(\cd ,z)\mapsto o^I$ given in (2)
  defines a linear isomorphism of linear spaces $I\binom{W^2}{U\ W^1}$
  and $\hom_{A_t(V)}$ $\big( A_t(U)\tensor_{A_t(V)}$ $W^1(0)$,
  $W^2(0)\big)$. 
\end{thm}

\pf
One can find similar proof in \cite{Li2} and \cite{Y}.
\qed

\begin{rem}
  In order to prove the above assertions rigorously, we have to
  introduce some cocycles for Lie superalgebra associated to an
  affinized SVOA. However, such  modifications can be cleared by a
  little attentions.
\end{rem}

\section{Some facts on SVOAs} 

In this section we give some notes on the basic property of SVOA and
its representations.

\subsection{Invariant bilinear form} 

Let $V$ be an SVOA and $M=\oplus_{n\in \hf\N} M(n)$ be its module.
We can find a natural $V$-module structure in the dual space
$M^*=\oplus_{n\in\hf \N} M(n)^*$.

\begin{df}\label{invariant_form}
  For the restricted dual space $M^*$ of a $\Z_2$-graded $V$-module
  $M$, we define the {\it adjoint vertex operators} $Y^*(a,z)$ by
  means of the linear map
  $$
  \begin{array}{l}
    V\ni a\mapsto Y^*(a,z)=\dsum_{n\in\Z} a^*_n z^{-n-1}\in \End (M^*)
    [[z,z^{-1}]]
  \end{array}
  $$
  determined by the condition
  \begin{equation} \label{invariant form}
    \lfm Y^*(a,z) f | v\rfm := \lfm f | Y(e^{zL_1} (\lambda
    z^{-2})^{L_0} \lambda^{-2L_0^2} a,z^{-1}) v\rfm
  \end{equation}
  for $a\in V$, $f\in M^*$ and $v\in M$, where $\lambda =e^{\pm \pii}$.    
\end{df}

\begin{rem}
  In the above definition, it seems that we have two definitions
  according to the choice of the square root of unity $\lambda$, 
  but each choice determines the same adjoint vertex operators
  since we have assumed that $V$ has $\hf\Z$-grading.
  So we may choose each root of unity.
\end{rem}

Similar to the case of VOAs, we have the following.

\begin{prop} (Theorem 5.2.1 in \cite{FHL})
  $(M^*, Y^*)$ is a $\Z_2$-graded $V$-module. Moreover, if each
  homogeneous space $M(n)$ of $M$ is finite dimensional, then
  $(M^*)^*\simeq M$.
\end{prop}

The proof is the same as that in \cite{FHL}. One can also show the
following.

\begin{prop}
  (1)\ (Proposition 5.3.6 in \cite{FHL})
  Assume that $V$ as $V$-module is isomorphic to $V^*$.
  Then the natural pairing $V\times V\to \C$ is automatically
  symmetric. 
  \vsb\\
  (2)\ (Theorem 3.1 in \cite{Li3})
  The space of invariant forms on $V$ is linearly isomorphic to 
  $\hom_\C \big( V_0/L_1V_1, \C \big)$. Therefore, the existence of
  invariant forms on $V$ is equivalent to the existence of invariant
  forms on sub VOA $V^0$ of $V$.
\end{prop}

We can also introduce the dual module for every $\Z_2$-twisted
$V$-module in similar way. 
However, we won't give the details here.

\subsection{$\Z_2$-conjugacy} 

Every SVOA $V=V^0\oplus V^1$ has a canonical involution $\sigma$ 
which is identical on $V^0$ and acts as $-1$ on $V^1$, so we can think 
of the $\sigma$-conjugation of $V$-modules.
That is, for a $V$-module $(M, Y_M)$, define another vertex operator 
$Y_M^\sigma$ by
$$
  Y_M^\sigma (a,z):= Y_M(\sigma a,z).
$$
Then $(M,Y_M^\sigma )$ is also a $V$-module and we will denote it
simply by $M^\sigma$.
However, $\sigma$-conjugation of a $\Z_2$-graded $V$-module is a
trivial concept because every $\sigma$-conjugate is isomorphic to
original $\Z_2$-graded $V$-module as one can easily see.
But, as the following proposition insists, the $\sigma$-conjugation of 
$\Z_2$-twisted $V$-modules is not a trivial concept.

\begin{prop}\label{Z_2-conj}
  Let $V=V^0\oplus V^1$ be a simple SVOA with $V^1\ne 0$ and
  $W$ be an irreducible $\Z_2$-twisted $V$-module.
  Then one of the followings holds:
  \vsb\\
  (1)\ $W$ and $W^\sigma$ are non-isomorphic $V$-modules if $W$ is
  irreducible as $V^0$-module.
  \vsb\\
  (2)\ If $W$ is not an irreducible $V^0$-module, then $W$ is
  completely reducible $V^0$-module and it has two irreducible
  components. Write $W=N^1\oplus N^2$, then we have $V^1\cd N^1=N^2$
  and $V^1\cd N^2=N^1$. Therefore, $W$ has a $\Z_2$-grading under the
  action of $V$ and the conjugate module $W^\sigma$ is isomorphic to
  $W$ as $V$-module.
\end{prop}

\pf
Assume that $W$ as $V^0$-module is irreducible and $W$ and $W^\sigma$
are isomorphic as $V$-modules.
Then we can find a $V$-isomorphism $f: W\to W^\sigma$.
Clearly, $f$ is also a $V^0$-isomorphism.
By definition, we have a linear map $\phi : W\to W^\sigma$ satisfying
$\phi Y(a,z)w=Y(\sigma a,z)\phi w$ for $a\in V$ and $w\in W$.
Then the linear map $\phi^{-1} f$ becomes a $V^0$-isomorphism of
$W$ so that we can find some non-zero scalar $\alpha\in \C$ 
such that $\phi^{-1} f=\alpha$ by Schur's lemma.
Let $a\in V^1$ and $w\in W$ be non-zero elements.
It is well-known that $Y(a,z)w\ne 0$ for simple SVOAs so that 
we get $0\ne f(Y(a,z)w)=\alpha \phi Y(a,z)w =-\alpha Y(a,z) \phi (w)
= -Y(a,z) f(w)=-f(Y(a,z)w)$, a contradiction.
Therefore $W$ and $W^\sigma$ are not isomorphic $V$-modules.
This proves (1) and so we may assume that $W$ is not an irreducible
$V^0$-module.
Let $N$ be a proper $V^0$-submodule of $W$.
Then the associativity implies $V^1\cd N$ is also $V^0$-submodule of
$W$ and we have $N\cap (V^1\cd N)=0$ since it is a proper
$V$-submodule of $W$.  
Therefore $W$ contains $N\oplus V^1\cd N$ and by the irreducibility we
obtain $W=N\oplus V^1\cd N$.  
This implies that any proper $V^0$-submodule $N$ of $W$ is irreducible 
$V^0$-module. Therefore $W$ has a $\Z_2$-grading under the action of
$V^0\oplus V^1$ and so the conjugate $W^\sigma$ is isomorphic to $W$
as $V$-module.
\qed

\section{Ising model SVOA} 

In this section we will give an explicit construction of the Ising
model SVOA $\ising\oplus\isingodd$ and its $\Z_2$-twisted modules
$\isingmod^\pm$. By calculating Zhu algebras explicitly, we shall
prove that $\isingmod^\pm$ are all irreducible $\Z_2$-twisted
modules for this SVOA.
This construction is well-known and the most of contents in Sec.{}
4.1-4.2 can be found in \cite{FRW} and \cite{FRW2}.

\subsection{Realization of Ising models} 

Here we consider a certain realizations of unitary highest weight
representations of Virasoro algebras of central charge $\fr{1}{2}$.
There are exactly three unitary representations $\ising$, $\isingodd$
and $\isingmod$, which are often called the {\it Ising models}. 
First two ones are realized as follow. 

Let $\mathcal{A}_\psi$ be the algebra generated by $\{ \psi_k | k\in 
\Z+\fr{1}{2}\}$ subject to the defining relations
$$
\begin{array}{ll}
  [\psi_m,\psi_n]_+:=\psi_m\psi_n+\psi_n\psi_m =\delta_{m+n,0},
  & m,n \in \Z +\fr{1}{2},
\end{array}
$$
and denote a subalgebra of $\mathcal{A}_\psi$ generated by $\{ 
\psi_k | k\in \Z +\fr{1}{2}, k>0\}$ by $\mathcal{A}_\psi^+$.
Let $\C \vac$ be a trivial $\mathcal{A}_\psi^+$-module.
Define a canonical induced $\mathcal{A}_\psi$-module $M$ by
$$
  M:=\mathrm{Ind}_{\mathcal{A}_\psi^+}^{\mathcal{A}_\psi} \C \vac
  = \mathcal{A}_\psi \tensor_{\mathcal{A}_\psi^+} \C \vac .
$$
As well-known, we can find Virasoro module structure in $M$.
Following \cite{KR}, set 
\begin{equation}\label{vir}
  L_n:=\fr{1}{2}\dsum_{k>-n/2}(n+2k)\psi_{-k}\psi_{n+k}, 
  \q (n\in\Z ).
\end{equation}
Then $\{ L_n | n\in\Z\}$ gives a representation of Virasoro algebra
of central charge $\fr{1}{2}$ on $M$ and as Virasoro module $M$
is decomposed as follow
$$
  M=\ising \oplus \isingodd ,
$$
where $L(c,h)$ denotes the irreducible highest weight 
module of central charge $c$ with highest weight $h$.
The component $\ising$ is generated by $\vac$, whereas $\isingodd$ is
generated by $\psi_{-\fr{1}{2}}\vac$ under the Virasoro and it is also
clear that the above decomposition coincides with the standard
$\Z_2$-grading decomposition, i.e. 
$$
  L\big( \tfr{1}{2},\tfr{i}{2}\big) 
  = \C \big\lfm\ \psi_{-n_1}\cds \psi_{-n_k}\vac\ \big|\ n_j>0,\ 
  n_1+\cds +n_k\in\Z+\tfr{i}{2}\ \big\rfm
$$
for $i=0,1$. Furthermore, one can introduce a symmetric contravariant 
Hermitian form $\lfm \cd | \cd \rfm$ on $M$ such that $\lfm \vac |
\vac \rfm =1$ and $\lfm \psi_n a | b\rfm = \lfm a|\psi_{-n} b\rfm$.
We also note that the above basis forms an orthonormal basis for $(M,
\lfm \cd | \cd \rfm )$.

Similarly, $\isingmod$ shall be realized as follow.
Let $\mathcal{A}_\phi$ be the other algebra
generated by $\{ \phi_n | n\in\Z\}$ whose defining relations are
$$
  [\phi_m,\phi_n]_+=\delta_{m+n,0},\q m,n\in\Z .
$$
Let $\mathcal{A}_\phi^+$ be a subalgebra of $\mathcal{A}_\phi$
generated by $\{ \phi_n | n>0\}$ and denote a trivial 1-dimensional 
$\mathcal{A_\phi^+}$-module by $\C v_0$. 
Then set $N=\mathrm{Ind}_{\mathcal{A}_\phi^+}^{\mathcal{A}_\phi}\C
v_0$ as we did previously. Again we can find Virasoro
representation in $N$. Set 
\begin{equation}\label{vir'}
  L_n':= \fr{1}{16} \delta_{n,0}
       + \fr{1}{2} \dsum_{k>-n/2} (n+2k)\phi_{-k}\phi_{n+k}, 
       \q  (n\in\Z ).
\end{equation}
Then $\{ L_n' | n\in\Z\}$ satisfies the Virasoro relation of central
charge $\fr{1}{2}$ on $N$. There are two distinct highest weight
vectors $v_0$ and $\phi_0v_0$ with highest weight $\fr{1}{16}$ in $N$
and as Virasoro module $N$ decomposes as follow (cf. \cite{KR}):
$$
  N=\isingmod\oplus\isingmod ,
$$
where one of $\isingmod$ is generated by $v_0$ and the other one is
generated by $\phi_0v_0$ under the Virasoro.
One can also introduce a symmetric contravariant Hermitian form $\lfm
\cd | \cd \rfm$ on $N$ such that $\lfm v_0 | v_0 \rfm = 1$, $\lfm v_0
| \phi_0 v_0\rfm = \lfm \phi_0 v_0 | v_0\rfm =0$ and $\lfm \phi_n a |
b\rfm = \lfm a | \phi_{-n}b\rfm$.

\subsection{SVOA structure on Ising models} 

We keep the same notation as previous.
By its construction, $M$ is generated by $\vac$ over
$\mathcal{A}_\psi$. Define the generating series
$$
  \psi (z):=\dsum_{n\in\Z}\psi_{n+\fr{1}{2}}z^{-n-1}.
$$
Since $[\psi (z),\psi (w)]_+ = z^{-1}\delta (\fr{w}{z})$, $\psi (z)$ 
is local with itself and it follows from the defining relations of
$\mathcal{A}_\psi$ that $\psi (z)$ satisfies $L_{-1}$-derivation
property $[L_{-1},\psi (z)]=\fr{d}{dz}\psi (z)$. 
Therefore we can consider a subalgebra of a local system on $M$
generated by $\psi (z)$ and $I(z)=\id_M$ (cf. \cite{Li1}).
By a direct calculation, one sees that 
$$
  \fr{1}{2}\psi (z)\circ_{-2}\psi (z)=\dsum_{n\in\Z} L_n z^{-n-2},
$$
where $\circ_n$ denotes the $n$-th normal product (cf. \cite{Li1}).
That is, the coefficients in the left hand side of the above equality 
coincide with \eqref{vir}. 
Recall that $\psi_{-n_1+\fr{1}{2}} \psi_{-n_2+\fr{1}{2}} \cds 
\psi_{-n_k+\fr{1}{2}}\vac$, $n_1>n_2>\cds >n_k>0$, $k\geq 0$ form a
basis of $M$. 
We shall define a vertex operator of each base $\psi_{-n_1+\fr{1}{2}}
\psi_{-n_2+\fr{1}{2}} \cds \psi_{-n_k+\fr{1}{2}}\vac$ on $M$.
For $k=0$ we set $Y(\vac,z):=\id_M$ and inductively
we define $Y(\psi_{-n+\fr{1}{2}}a,z):=\psi (z)\circ_n Y(a,z)$.
Then by the theory of the local system, we have the following well-known
statement.

\begin{thm} (Theorem 2 in \cite{FRW2})\
  By the above definition, $(\ising\oplus\isingodd, Y(\cd,z),$ 
  $\vac,\fr{1}{2}\psi_{-\fr{3}{2}}\psi_{-\fr{1}{2}}\vac )$ has a simple 
  (and unique) SVOA structure with $\ising$ even part and $\isingodd$
  odd part. 
\end{thm}

Note that the invariant bilinear form defined by \eqref{invariant form} 
coincides with the contravariant bilinear form on $M$ since
$\psi_n^*=\psi_{-n}$ for all $n\in \Z +\hf$.
Next, we consider $\ising\oplus\isingodd$-module structure in
$\isingmod$. Recall that in $N$ there are two highest weight vectors 
$v_0$ and $\phi_0v_0$ and each of them generates $\isingmod$ under
the Virasoro algebra.
Since they are not stable under $\phi_0$, we have to change the
highest weight vectors to suitable ones.  
We define $v^\pm_{\fr{1}{16}}:=\phi_0 \vac \pm
\fr{1}{\sqrt{2}} \vac$ to have $\phi_0 \cd v^\pm_{\fr{1}{16}} = \pm
\fr{1}{\sqrt{2}} v^\pm_{\fr{1}{16}}$. 
Then both $v^+_{\fr{1}{16}}$ and $v^-_{\fr{1}{16}}$ are highest weight
vectors and each of them generates $\isingmod$ under the Virasoro.
Denote the Virasoro modules generated by $v^\pm_{\fr{1}{16}}$ by
$\isingmod^\pm$, respectively. 
Then we have $N = \isingmod^+ \oplus \isingmod^-$.  
Note that $\isingmod^+$ and $\isingmod^-$ are isomorphic as Virasoro
modules but they are not isomorphic to each other as
$\mathcal{A}_\phi$-modules.
Consider the generating series
$$
  \phi (z):= \dsum_{n\in\Z}\phi_n z^{-n-\fr{1}{2}}.
$$
By direct calculations one can show that $\phi (z)$ is local with
itself and the derivation property $[L_{-1}',\phi (z)] = \fr{d}{dz}
\phi (z)$ holds.  
We realize Virasoro operator \eqref{vir'} by using $\phi (z)$.
Since the powers of $z$ in $\phi (z)$ lie in $\Z +\fr{1}{2}$, we have
to use the twisted normal product. 
Define a generating series $L(z)$ of operators
on $N$ by
$$
\begin{array}{l}
     L(z_2)
     := \dfr{1}{2}\res_{z_0} \res_{z_1} z_0^{-2} \l( 
      \dfr{z_1-z_0}{z_2} \r)^{\fr{1}{2}} 
      \vsb\\
     \times \l\{ z_0^{-1}\delta\l(\dfr{z_1-z_2}{z_0}\r) \phi (z_1)  
      \phi (z_2) + z_0^{-1}\delta\l(\dfr{-z_2+z_1}{z_0}\r) \phi (z_2) 
      \phi (z_1)\r\} .
\end{array}
$$
The following is a consequence of a direct computation.

\begin{lem}
  \q $L(z)=\dsum_{n\in\Z} L_n'z^{-n-2}$, where $L_n'$ is defined by
  \eqref{vir'}.
\end{lem}

Thanks to the above lemma, we can find a $\Z_2$-twisted $\ising \oplus
\isingodd$-module structure in $\isingmod^\pm$.
To show it, we associate a vertex operator on $N$ for every element of 
$\ising \oplus \isingodd$ and then we prove that our association gives
a homomorphism of vertex superalgebras. 
Set $Y_N(\vac ,z):=\id_N$ and define inductively a vertex
operator of $\psi_{-n+\fr{1}{2}}a$ on $N$ by
$$
\begin{array}{l}
   Y_N(\psi_{-n+\fr{1}{2}}a,z)
   := \dfr{1}{2}\res_{z_0} \res_{z_1} z_0^{-n} \l( 
   \dfr{z_1-z_0}{z_2} \r)^{\fr{1}{2}} 
   \vsb\\
   \times \l\{ z_0^{-1}\delta\l(\dfr{z_1-z_2}{z_0}\r) \phi (z_1)  
   Y_N(a,z_2) -(-1)^{q(a)} z_0^{-1} \delta\l(\dfr{-z_2+z_1}{z_0}\r) 
   Y_N(a,z_2) \phi (z_1)\r\} ,
\end{array}
$$
where $a=\psi_{-n_1+\fr{1}{2}}\cds \psi_{-n_k+\fr{1}{2}}\vac$, 
$n>n_1>\cds >n_k>0$ and extend linearly on $\ising \oplus\isingodd$.
Let $\mathfrak{A}$ be a $\Z_2$-twisted local system on $N$ in which
$\phi (z)$ contained. 
It is shown in \cite{Li1} that $\mathfrak{A}$ is a vertex superalgebra
under the $\Z_2$-twisted normal product (Theorem 3.14 in \cite{Li1}).
See \cite{Li1} for the detailed description of the twisted
normal product. 

\begin{lem}\label{SVA_hom}
  The linear mapping $a\in \ising\oplus\isingodd \mapsto Y_N(a,z)\in
  \mathfrak{A}$ defined above gives an vertex superalgebra
  homomorphism. 
\end{lem}

\pf
We should show that $Y_N(a_mb,z)=Y_N(a,z)\circ_m Y_N(b,z)$ for any
$a,b\in \ising\oplus\isingodd$ and $m\in\Z$, where $\circ_m$ denotes
the $\Z_2$-twisted $m$-th normal product in $\mathfrak{A}$. 
We may assume that $a=\psi_{-n_1+\fr{1}{2}}\cds 
\psi_{-n_k+\fr{1}{2}}\vac$, $n_1>\cds >n_k>0$. We proceed by induction 
on $k$. The case $k=0$ is trivial and the case $k=1$ is just the
definition. Assume that $k\geq 1$ and $Y_N(a_mb,z)= Y_N(a,z)\circ_m
Y_N(b,z)$ holds for arbitrary $b\in \ising \oplus \isingodd$ and
$m\in\Z$. Take any $n>n_1$. Then we have
$$
\begin{array}{l}
  Y_N\l( (\psi_{-n+\fr{1}{2}}a)_mb,z\r) 
  \vsb\\
  = \dsum_{i=0}^\infty (-1)^i \dbinom{-n}{i} Y_N \l( 
    \psi_{-n-i+\fr{1}{2}} a_{m+i} b  - (-1)^{q(a)-n} a_{-n+m-i}
    \psi_{i+\fr{1}{2}}b,z\r)
  \vsb\\
  = \dsum_{i=0}^\infty (-1)^i \dbinom{-n}{i} \Big\{ \phi (z)
    \circ_{-n-i} Y_N( a_{m+i} b,z) 
  \vsb\\
  \hspace{3cm}
    -(-1)^{q(a)-n} Y_N(a,z) 
    \circ_{-n+m-i} Y_N( \psi_{i+\fr{1}{2}} b,z) \Big\}
  \vsb\\
   = \dsum_{i=0}^\infty (-1)^i \dbinom{-n}{i} \Big\{ \phi (z)
    \circ_{-n-i} \Big( Y_N(a,z)\circ_{m+i} Y_N(b,z) \Big)
   \vsb\\
   \hspace{3cm}  
     -(-1)^{q(a)-n} Y_N(a,z) \circ_{-n+m-i} \Big( \phi (z)
     \circ_i Y_N( b,z)\Big) \Big\}
  \vsb\\
  = \Big( \phi (z)\circ_{-n} Y_N(a,z)\Big) \circ_m Y_N(b,z)
  \qq \text{(associativity\ in}\ \mathfrak{A})
  \vsb\\
  = Y_N\l( \psi_{-n+\fr{1}{2}}a,z \r) \circ_m Y_N(b,z).
\end{array}
$$
Therefore by induction the mapping $a\mapsto Y_N(a,z)$ defines a
vertex superalgebra homomorphism.
\qed
\vsv\\
Ler $V$ be an arbitrary SVOA.
By Proposition 3.17 in \cite{Li1},
giving a $\Z_2$-twisted $V$-module structure on $N$ is equivalent to
giving a vertex superalgebra homomorphism from $V$ to some local
system of $\Z_2$-twisted vertex operators on $N$.
Since both $\isingmod^+$ and $\isingmod^-$ are stable under the
action $Y_N(\cd\, ,z)$, we arrive at the following conclusion.

\begin{thm}\label{1/16pm} (Theorem 2 in \cite{FRW2})\  
  $\mathcal{A}_\phi$-modules $\isingmod^+$ and $\isingmod^-$ are
  non-isomorphic irreducible $\Z_2$-twisted
  $\ising\oplus\isingodd$-modules. 
\end{thm}

Note that the canonical contravariant bilinear form on $N$ can be
naturally extended to the invariant bilinear form for the module
vertex operators defined by \eqref{invariant form}.
Therefore, by the natural symmetries of the fusion rules, we can
explicitly calculate all intertwining operations of any type for the 
Ising models. 

\begin{rem}
  The $\Z_2$-twisted module vertex operators $Y_N(\cd\,
  ,z)|_{\isingmod^\pm}$ defined above give intertwining
  operators of type $L(\fr{1}{2},h)\times \isingmod^\pm\to
  \isingmod^\pm$ for $h=0,\fr{1}{2}$, respectively. Therefore, our
  construction gives another proof of Proposition 4.1-4.2 in
  \cite{M2}. 
\end{rem}

\begin{rem}\label{top}
  The $\Z_2$-twisted $\ising \oplus \isingodd$-modules 
  $\isingmod^\pm$ can be characterized as follow.   
  On the top level of $\isingmod^+$ $(=\C v^+_{\fr{1}{16}})$,
  a highest weight vector $\psi_{-\fr{1}{2}}\vac$ in $\isingodd$
  acts as $\fr{1}{\sqrt{2}}$, whereas on the top level of
  $\isingmod^-$ it acts as $-\fr{1}{\sqrt{2}}$.
  This observation will be used to determine the $\Z_2$-twisted 
  Zhu algebra associated to the Ising model SVOA in the next subsection. 
\end{rem}

\subsection{Zhu algebras for Ising model SVOA}

\gs{In the following context, we will use the highest weight 
to denote the Virasoro module itself to simplify the notation}.
For example, we shall denote $[0]$ for $\ising$, $[\fr{1}{2}]$ for
$\isingodd$, and so on.

Recall the following theorem in \cite{KW}.

\begin{thm}\ (Theorem 4.1 in \cite{KW})\q \label{KW}
  $A\l( [0]\oplus [\fr{1}{2}]\r) \simeq \C$.
  Therefore, SVOA $[0]\oplus [\fr{1}{2}]$ has a unique irreducible
  $\Z_2$-graded representation, namely $[0]\oplus [\fr{1}{2}]$ itself.
\end{thm}

It is obvious that the canonical involution $\sigma$ is the only
non-trivial automorphism on $[0] \oplus [\fr{1}{2}]$. 
Therefore the $\Z_2$-twisted representation is the only 
non-trivial twisted representation of $[0] \oplus [\fr{1}{2}]$.
Let's determine the $\Z_2$-twisted Zhu algebra $A_t( [0] \oplus
[\fr{1}{2}] )$. Denote by $\iota$ a canonical surjection from
$[0] \oplus [\fr{1}{2}]$ onto $A_t( [0] \oplus [\fr{1}{2}])$.

\begin{lem}\label{lemma} \ 
  $\iota \l( [0] \r) =\C \iota (\vac )$ and $\iota\l( [\fr{1}{2}] \r)
  = \C \iota (\psi_{-\fr{1}{2}}\vac )$.
\end{lem}

\pf
Again we proceed by induction on the length $k$ of a canonical
basis  
$$
  \psi_{-n_1+\fr{1}{2}} \cds \psi_{-n_k+\fr{1}{2}}\vac 
$$
with $n_1>\cds>n_k>0$. 
The case $k=0$ is obvious. For $k=1$, by Lemma
\ref{affinezhu} we have
$$
  \res_z Y(a,z)b\dfr{(1+z)^{\wt (a)}}{z^{2+m}}\in O_t \l( [0] 
  \oplus [\tfr{1}{2}] \r) 
$$
for any $a,b\in [0] \oplus [\fr{1}{2}]$ and arbitrary $m\geq 0$.
Putting $a=\psi_{-\fr{1}{2}}\vac$, $b=\vac$ and by induction on
$n$, we obtain $\psi_{-n+\fr{1}{2}}\vac \in O_t\l( [0]\oplus
[\fr{1}{2}] \r)$ for $n\geq 2$. Thus the case $k=1$ is correct.
Assume that our assertion is correct for all
$\psi_{-n_1+\fr{1}{2}}\cds \psi_{-n_s+\fr{1}{2}}\vac$ with $s\leq k$. 
Take any $n>n_1$, $n>1$. Then 
\begin{equation}\label{calc}
\begin{array}{l}
  O_t\l( [0]\oplus [\tfr{1}{2}] \r)\ni \res_z\psi (z)
    \dfr{(1+z)^{\fr{1}{2}}}{z^{2+(n-2)}} \psi_{-n_1+\fr{1}{2}}\cds 
    \psi_{-n_s+\fr{1}{2}}\vac 
  \vsb\\
  = \psi_{-n+\fr{1}{2}} \psi_{-n_1+\fr{1}{2}} \cds
    \psi_{-n_k+\fr{1}{2}} \vac 
  + \dsum_{i=1}^\infty \dbinom{\fr{1}{2}}{i} \psi_{-n+i+\fr{1}{2}}
    \psi_{-n_1+\fr{1}{2}} \cds \psi_{-n_k+\fr{1}{2}} \vac .
\end{array}
\end{equation}
It follows from the inductive assumption that the second term in the
lower hand side of \eqref{calc} can be written as the desired forms.
Hence the first term in the lower hand side also has the desired form
since each term in the lower hand side of \eqref{calc} shares the same
parity. So our assertion holds for $k+1$. Thus we have the desired
result. 
\qed

The following theorem together with Theorem \ref{KW} completes the
classification of all irreducible representations of the Ising model
SVOA.

\begin{thm}
  There exists exactly two non-isomorphic irreducible $\Z_2$-twisted 
  $[0]\oplus [\fr{1}{2}]$-modules, namely, $[\fr{1}{16}]^+$ and
  $[\fr{1}{16}]^-$.
\end{thm}

\pf
We have already shown that both $[\fr{1}{16}]^+$ and $[\fr{1}{16}]^-$
are irreducible $\Z_2$-twisted $[0]\oplus [\fr{1}{2}]$-modules in
Theorem \ref{1/16pm}.
By Theorem \ref{1:1}, we should determine all irreducible
$A_t([0]\oplus [\fr{1}{2}])$-modules.
By Lemma \ref{lemma}, we know that $A_t([0]\oplus [\fr{1}{2}])$ is
generated by $\iota (\vac )$ and $\iota (\psi_{-\fr{1}{2}}\vac)$.
After a short computation along with the definition, we get the
followings. 
$$
\begin{array}{ll}
  (\psi_{-\fr{1}{2}}\vac )\circ_t (\psi_{-\fr{1}{2}}\vac )
  = 2\l(\w -\dfr{1}{16}\vac \r) ,
  &
  (\psi_{-\fr{1}{2}}\vac )*_t (\psi_{-\fr{1}{2}}\vac ) =
  \dfr{1}{2}\vac .
\end{array}
$$
Hence we can deduce that 
$\w_1=L_0$ acts as $\fr{1}{16}$ on the top level of every irreducible 
$\Z_2$-twisted $[0]\oplus [\fr{1}{2}]$-module and
$A_t\l([0]\oplus[\fr{1}{2}]\r)$ is a homomorphic image of a ring $\C
[x]/\lfm x^2-\fr{1}{2}\rfm$. But by Remark \ref{top},
the top levels $\C v^\pm_{\fr{1}{16}}$ of $[\fr{1}{16}]^\pm$ are
non-isomorphic irreducible $A_t([0]\oplus [\fr{1}{2}])$-modules on
which $\iota (\psi_{-\fr{1}{2}}\vac )$ acts as $\pm \fr{1}{\sqrt{2}}$, 
respectively, so that $A_t([0]\oplus [\fr{1}{2}])$ must be isomorphic
to $\C [x]/\lfm x^2-\fr{1}{2}\rfm$. 
Therefore, $[\fr{1}{16}]^+$ and $[\fr{1}{16}]^-$ are only
non-isomorphic irreducible $\Z_2$-twisted $[0]\oplus
[\fr{1}{2}]$-modules.
\qed

\begin{rem}\label{scalar}
  Since $[\fr{1}{16}]$ as $[0]$-module is irreducible, we know
  that $([\fr{1}{16}]^\pm)^\sigma$ is not isomorphic to
  $[\fr{1}{16}]^\pm$ as $[0]\oplus [\hf]$-module, respectively,  by
  Proposition \ref{Z_2-conj}.  
  This implies that if we denote the module vertex operator of $[\hf
  ]$ on $[\fr{1}{16}]^+$ by $I(\cd ,z)$, then the module vertex
  operator of $[\hf ]$ on $[\fr{1}{16}]^-$ is given by $-I(\cd ,z)$.
  Since $I(\cd ,z)$ and $-I(\cd ,z)$ are the only scalar multiples of
  $I(\cd ,z)$ which satisfy the $\Z_2$-twisted Jacobi identity,
  we can also conclude that $[\fr{1}{16}]^\pm$ are all non-isomorphic
  $\Z_2$-twisted $V$-modules.
\end{rem}

\subsection{Application} 

In this subsection we consider an applications of the Ising models.
One of the merit of using the Ising model is that it defines an
automorphism, so-called ``Miyamoto involution'' of VOA.
Concerning to this involution, we will find some $\Z_2$-twisted
representations of SVOAs.

In the following, we will assume that $V=\oplus_{n=0}^\infty V_n$ is 
a simple VOA whose weight zero space $V_0$ is spanned by $\vac$, i.e.
$V_0=\C \vac$. 
Following \cite{M1}, we will call an element $e\in V_2$ whose vertex
operator $Y(e,z)=\sum_{n\in \Z} L^e_n z^{-n-2}$ generates a copy of
the Virasoro algebra of central charge $c_e$ a {\it conformal vector}
with c.c. $c_e$. One can find the following lemma in \cite{M1}.

\begin{lem}
  An element $e\in V_2$ is a conformal vector with c.c. $c_e$ if and
  only if it satisfies $e_1e=2e$, $e_2e=0$ and $e_3e=\fr{1}{2}c_e
  \vac$.
\end{lem}

A conformal vector $e$ is called {\it rational} if the sub VOA $\lfm e 
\rfm$ generated by $e$ becomes a rational Virasoro VOA.
In this paper, we are especially interested in a rational conformal
vector $e$ with c.c. $\hf$.
One way to determine whether a conformal vector $e$ with c.c. $\hf$ is 
rational or not is to check $64(L^e_{-2})^3 \vac +93(L^e_{-3})^2\vac
-264 L^e_{-4}L^e_{-2} \vac -108 L^e_{-6}\vac =0$ or not
(cf. \cite{DMZ}). 
Assume that $e$ is a rational conformal vector with c.c. $\hf$ in $V$.
Seen $V$ as $\lfm e\rfm \simeq \ising$-module, we can decompose $V$ as
follow.
$$
  V=V_e(0)\oplus V_e(\tfr{1}{2} )\oplus V_e(\tfr{1}{16}),
$$
where $V_e(h)$ is the sum of all irreducible $\lfm e\rfm$-modules
which are isomorphic to $L(\hf, h)$ for $h=0,\hf,\fr{1}{16}$ as
Virasoro modules.
For $h=0,\hf,\fr{1}{16}$, let $T_e(h)$ be the space $\{ v\in V |
L^e_0v =hv\}$ of the highest weight vectors.
Then $T_e(h)$ is also the space of the multiplicity of irreducible
components in $V$ and so we get a linear isomorphism $V_e(h)\simeq
L(\hf,h)\tensor_\C T_e(h)$.
Since $[e_1,\w_1]=0$, $T_e(h)$ admits a $L_0$-weight space
decomposition $T_e(h)=\oplus_{n\in \N} T_e(h)_n$, where
$T_e(h)_n=V_n\cap T_e(h)$. 
By a little calculation, one can show that $\w_2 e=0$. So by \cite{FZ}
we know that $f:=\w -e$ is also a conformal vector with
c.c. $\mathrm{rank}(V)-\hf$ whose vertex operators are commutative with 
that of $e$.
Therefore we have a decomposition $\w =e+f$ of the Virasoro vector of
$V$ into two mutually commutative conformal vectors.
It is also shown in \cite{FZ} that if we set the commutants  $U^e:=
\ker_V L^f_{-1}$ and $U^f:= \ker_V L^e_{-1}$, then $(U^e, Y_{|U^e},
\vac_V, e)$ and $(U^f, Y_{|U^f}, \vac_V, f)$ forms mutually
commutative sub VOAs of $V$, where we set $L^e_n=e_{n+1}$ and
$L^f_n=f_{n+1}$. 
In general, $\ker_V L^e_{-1} \subseteq \ker_V L^e_0$ and
$\ker_V L^f_{-1} \subseteq \ker_V L^f_0$, but in our case the equality
holds for both cases and we have $U^e=\lfm e\rfm$ and $U^f=T_e(0)$.
Furthermore, we can show the following.

\begin{prop}
  $(T_e(0), Y_{|T_e(0)}, \vac_V, f=\w -e)$ is the maximal simple sub
  VOA of $V$ whose vertex operators are commutative with that of $\lfm
  e\rfm$ element-wisely in $V$.
  Moreover, $V_e(0)=\ising\tensor T_e(0)$ is a simple sub VOA of $V$
  whose Virasoro vector is the same as that of $V$.
\end{prop}

\pf
We should only show the simplicity of $T_e(0)$.
On $V$, the linear map $\tau_e$ which is identical on
$V_e(0)\oplus V_e(\hf )$ and acts $-1$ on $V_e(\fr{1}{16})$ defines an 
involution of $V$, which is known as ``Miyamoto automorphism'' in
\cite{M1}. 
On the $\tau_e$-invariants $V_e(0)\oplus V_e(\hf )$, the linear map
$\sigma_e$ which is identical on $V_e(0)$ and acts as $-1$ on $V_e(\hf 
)$ also defines an involution \cite{M1}.
Therefore the $\sigma_e$-invariants $V_e(0)=\ising \tensor T_e(0)$ is
a simple sub VOA of $V$ since $V$ is simple.
Then the irreducibility of $\ising$ implies $T_e(0)$ is also simple.
\qed

By this proposition, we may view $V$ as a $[0]\tensor T_e(0)$-module. 
It is not difficult to show that $T_e(h)$, $h=0,\hf,\fr{1}{16}$ are
$T_e(0)$-modules. Since the $\tau_e$-invariants $[0]\tensor
T_e(0)\oplus [\hf ]\tensor T_e(\hf )$ is a simple sub VOA 
of $V$, $T_e(\hf )$ is an irreducible $T_e(0)$-module.
On the other hand, $[0]\oplus [\hf ]$ has an SVOA structure. 
So it is natural for us to expect that $T_e(0)\oplus T_e(\hf )$ also
has an SVOA structure. 
Namely, we expect that the decomposition $V^{\tau_e}=[0]\tensor
T_e(0)\oplus [\hf ]\tensor T_e(\hf )$ describes not only the tensor
product of vector spaces but also that of SVOAs.
As a generalization of Proposition 4.9 in \cite{M1}, the following
theorem holds.

\begin{thm}\label{cutting SVOA}
  Assume that $V_e(\hf )\ne 0$. Then there exists a simple SVOA
  structure on $T_e(0)\oplus T_e(\hf )$ 
  such that the even part of a tensor product of SVOAs $[0]\oplus [\hf 
  ]$ and $T_e(0)\oplus T_e(\hf )$ is isomorphic to $V_e(0)\oplus
  V_e(\hf )$ as VOAs.
\end{thm}

\pf
We shall define vertex operators on an abstract linear space
$T_e(0)\oplus T_e(\hf )$.
First, we show the existence of the
intertwining operators.
Let $\{ a^\gamma | \gamma \in \Gamma\}\subset T_e(0)$ and $\{ 
u^\lambda | \lambda \in \Lambda\}\subset T_e(\hf )$ be sets of basis
elements of $T_e(0)$ and $T_e(\hf )$ consisting of $L^f_0$-eigen
vectors, respectively.
Using these bases, we obtain decompositions of $V_e(0)$ and $V_e(\hf )$
as $[0]$-modules as follows.
$$
  V_e(0)=\bigoplus_{\gamma\in\Gamma} [0]\tensor a^\gamma,\q
  V_e(\thf )=\bigoplus_{\lambda\in\Lambda} [\thf ]\tensor u^\lambda .
$$
Let $\pi_\gamma$ be a projection map $V_e(0)\to [0]\tensor a^\gamma$.
We define intertwining operators $I_{\lambda\mu}^\gamma (\cd ,z)$ of
type $[\hf ]\tensor u^\lambda \times [\hf ]\tensor u^\mu \to
[0]\tensor a^\gamma$ as follow.
$$
  I_{\lambda\mu}^\gamma (x,z)y
  := \pi_\gamma Y(x\tensor u^\lambda ,z)
  y\tensor u^\mu z^{-\abs{\gamma}+\abs{\lambda}+\abs{\mu}},
$$
where $x,y\in [\hf ]$ and $\abs{\gamma}$, $\abs{\lambda}$, $\abs{\mu}$
denote the $L^f_0$-weights of $a^\gamma$, $u^\lambda$, $u^\mu$,
respectively.
One can show that $I_{\lambda\mu}^\gamma (\cd ,z)$ satisfy
$L^e_{-1}$-derivation, the condition of finiteness of negative powers
and Jacobi identity for intertwining operators so that
$I_{\lambda\mu}^\gamma(\cd ,z)$ are intertwining operators of type 
$[\hf ]\times [\hf ]\to [0]$.
Since the space of intertwining operators of such a type is one
dimensional, there exist suitable scalars $c_{\lambda\mu}^\gamma\in
\C$ such that $I_{\lambda\mu}^\gamma (\cd ,z)$ are a multiple of
the vertex operator $Y_{[\hf ]\times [\hf ]}(\cd ,z)$ by
$c_{\lambda\mu}^\gamma$, which was  constructed in Sec. 4.2.
Therefore, the vertex operator of $x\tensor u^\lambda$ on $[\hf
]\tensor T_e(\hf )$ can be written as follow.
$$
  Y_V(x\tensor u^\lambda ,z) y\tensor u^\mu 
  = Y_{[\hf ]\times [\hf ]}(x,z)y \tensor \dsum_{\gamma\in\Gamma}
    c_{\lambda\mu}^\gamma a^\gamma
    z^{\abs{\gamma}-\abs{\lambda}-\abs{\mu}} .
$$
Thus, by setting $J(u^\lambda ,z)u^\mu := \dsum_{\gamma\in \Gamma} 
c_{\lambda\mu}^\gamma a^\gamma z^{\abs{\gamma} -\abs{\lambda}
-\abs{\mu}}$, we obtain a decomposition 
$$
  Y_V(x\tensor u^\lambda ,z)y\tensor u^\mu =
  Y_{[\hf ]\times [\hf ]}(x,z)y \tensor J(u^\lambda ,z)u^\mu
$$ 
for $x\tensor u^\lambda, y\tensor u^\mu\in [\hf ]\tensor T_e(\hf )$.
We claim that $J(\cd ,z)$ is an intertwining operator of type $T_e(\hf 
)\times T_e(\hf )\to T_e(0)$.
It is not difficult to show that $J(\cd ,z)$ satisfies the condition
of finiteness of negative powers and $L^f_{-1}$-derivation.
So we should show the commutativity and associativity of $J(\cd ,z)$.
Let $a\in T_e(0)$, $u,v\in T_e(\hf )$ be arbitrary elements.
Take a sufficiently large $N\in \N$. Then the commutativity of vertex
operators on $V$ leads
$$
\begin{array}{l}
  (z_1-z_2)^N Y_V(\vac \tensor a,z_1) 
    Y_V(\psi_{-\hf}\vac\tensor u,z_2) \psi_{-\hf}\vac \tensor v
  \vsb\\
  = (z_1-z_2)^N Y_V(\psi_{-\hf}\vac \tensor u,z_1) 
    Y_V(\vac \tensor a,z_2) \psi_{-\hf}\vac \tensor v .
\end{array}
$$
Rewriting the above equality yields
$$
\begin{array}{l}
  (z_1-z_2)^N Y_{[\hf ]\times [\hf ]}(\psi_{-\hf}\vac ,z_2)
    \psi_{-\hf}\vac \tensor Y_{T_e(0)}(a,z_1) J(u,z_2)v
  \vsb\\
  = (z_1-z_2)^N Y_{[\hf ]\times [\hf ]}(\psi_{-\hf}\vac,z_2)
    \psi_{-\hf} \vac \tensor J(u,z_2)Y_{T_e(0)}(a,z_1)v.
\end{array}
$$
By comparing the coefficients of $(\psi_{-\hf}\vac )_0\psi_{-\hf}\vac
= \vac$, we get the commutativity:
$$
  (z_1-z_2)^N Y_{T_e(0)}(a,z_1) J(u,z_2)v
  = (z_1-z_2)^N J(u,z_2) Y_{T_e(0)}(a,z_1)v.
$$
Similarly, by considering some coefficients of $Y_V(Y_V(\vac \tensor
a,z_0)\psi_{-\hf}\vac \tensor u,z_2)\psi_{-\hf}\vac \tensor v$ in $V$, 
we can obtain the associativity.
Hence, $J(\cd ,z)$ is the intertwining operator of the desired type.

Using $J(\cd ,z)$, we introduce the vertex operations on $T_e(0)\oplus 
T_e(\hf )$.
Since we already know the action of $T_e(0)$ on $T_e(0)\oplus T_e(\hf
)$, we should define the action of $T_e(\hf )$ on $T_e(0)\oplus
T_e(\hf )$. For $a\in T_e(0)$, $u,v\in T_e(\hf )$, we set
$$
  Y_{T_e(\hf )}(u,z)a:= e^{zL^f_{-1}} Y_{T_e(0)}(a,-z)u,\q
  Y_{T_e(\hf )}(u,z)v:= J(u,z)v.
$$
Since the above vertex operators are intertwining operators of type
$T_e(\hf )\times T_e(0)\to T_e(\hf )$ and $T_e(\hf )\times T_e(\hf )\to
T_e(0)$, respectively, we only need to show the mutually commutativity of
vertex operators.
It follows from the definition that the vertex operator $Y_V(a\tensor
b ,z)$ in $V^{\tau_e}$ can be written as $Y_{[0]\oplus [\hf
]}(a,z)\tensor Y_{T_e(0)\oplus T_e(\hf )}(b,z)$ for $a\in [0]\oplus
[\hf ]$ and $b\in T_e(0)\oplus T_e(\hf )$. Hence, by comparing
suitable coefficients as we did previously, we can deduce the mutually
commutativity of the vertex operators on $T_e(0)\oplus T_e(\hf )$
since the vertex operators on $V$ and that on $[0]\oplus [\hf ]$
satisfy the mutually commutativity. 
Therefore, by our definition, $(T_e(0)\oplus T_e(\hf ), Y(\cd
,z),\vac, f)$ becomes a simple SVOA.
The rest of assertion is now clear.
\qed

Since $\tau_e^2=1$ on $V$, the space $V_e(\fr{1}{16})$ is irreducible
$V^{\tau_e}$-module. As $(V^{\tau_e})^{\sigma_e}=[0]\tensor
T_e(0)$-module, $V_e(\fr{1}{16})$ can be written as
$[\fr{1}{16}]\tensor T_e(\fr{1}{16})$. 
However, it is not clear that $T_e(\fr{1}{16})$ is irreducible under
$T_e(0)$ in general.
But we can insure that $T_e(\fr{1}{16})$ is irreducible under
$T_e(0)\oplus T_e(\hf )$.

\begin{thm} \label{Z_2 piece}
  Assume that $V_e(\fr{1}{16})\ne 0$. 
  Then $T_e(\fr{1}{16})$ has an irreducible $\Z_2$-twisted
  $T_e(0)\oplus T_e(\hf )$-module structure such that 
  $V_e(\fr{1}{16})$ is isomorphic to a tensor product of
  irreducible $\Z_2$-twisted $[0]\oplus [\hf ]$-module
  $[\fr{1}{16}]^+$ and irreducible $T_e(0)\oplus T_e(\hf )$-module
  $T_e(\fr{1}{16})$. 
\end{thm}

\pf
Since we know how $[0]\oplus [\hf ]$ acts on $[\fr{1}{16}]^\pm$,
we can use the same strategy of the proof of Theorem \ref{cutting SVOA}.
Since $V_e(\fr{1}{16})=[\fr{1}{16}]\tensor T_e(\fr{1}{16})$ is an
irreducible $V^{\tau_e}=[0]\tensor T_e(0)\oplus [\hf ]\tensor T_e(\hf
)$-module, it is not hard to see that the vertex operator of 
$a\tensor b\in [h]\tensor T_e(h)$ on $[\fr{1}{16}]\tensor
T_e(\fr{1}{16})$ can be written as $Y_{[h]\times [\fr{1}{16}]} (a,z)
\tensor Y_{T_e(h) \times T_e(\fr{1}{16})} (b,z)$, where
$Y_{[h] \times [\fr{1}{16}]} (\cd ,z)$ is an intertwining operator of
type $[h] \times [\fr{1}{16}]\to [\fr{1}{16}]$ and $Y_{T_e(h)\times
T_e(\fr{1}{16})}(\cd ,z)$ is an intertwining operator of type
$T_e(h)\times T_e(\fr{1}{16})\to T_e(\fr{1}{16})$ for $h=0,\hf$,
respectively. 
As mentioned in Remark \ref{scalar}, by managing salar multiplications
we may take $Y_{[h]\times [\fr{1}{16}]}(\cd ,z)$ to be the
$\Z_2$-twisted module vertex operators of $[0]\oplus [\hf ]$ on
$[\fr{1}{16}]^+$ for $h=0,\hf$.
Then the Jacobi identity on $V$ and the $\Z_2$-twisted Jacobi identity
on $[\fr{1}{16}]^+$ implies that the intertwining operators 
$Y_{T_e(h)\times T_e(\fr{1}{16})}(\cd ,z)$ on $T_e(\fr{1}{16})$ also
satisfy the $\Z_2$-twisted Jacobi identity for $h=0,\hf$.
Thus, $T_e(\fr{1}{16})$ is an irreducible $\Z_2$-twisted $T_e(0)\oplus
T_e(\hf )$-module under the module vertex operators $Y_{T_e(h)\times
T_e(\fr{1}{16})}(\cd ,z)$ for $h=0,\hf$.
\qed

\section{The Babymonster SVOA $\VB$} 

In this section, we apply our result to the moonshine vertex operator
algebra, which is of course the most important example of
holomorphic VOAs.
In this paper we treat the moonshine vertex operator algebra
$V^\nat_\R$ over \gs{the real number field} constructed in \cite{M3}. 
As well-known, the full automorphism group of $V^\nat_\R$ is the
Monster sporadic simple group $\M$ (cf. \cite{M3}).
In $V^\nat_\R$, there are many rational conformal vectors with
c.c. $\hf$. 
It is shown in \cite{M1} that each rational conformal vector $e\in
V^\nat_\R$ defines an element $\tau_e$ of 2A conjugacy class of $\M$.
It is also shown in \cite{C} that this correspondence is one-to-one. 
In $\M$, we can find many sporadic simple groups.
Let $e\in V^\nat_\R$ be a rational conformal vector with c.c. $\hf$.
Then the centralizer $C_{\M}(\tau_e)$ is isomorphic to an
extension $\lfm \tau_e \rfm \cd \B$ of the Baby monster sporadic
simple group $\B$. 
So $\B$ acts on the $\tau_e$-invariants of $V^\nat_\R$ as an
automorphism group of a VOA.  
Motivated by this fact, H\"{o}hn studied the $\tau_e$-invariants of
$V^\nat_\R$ and found the Babymonster SVOA $\VB$ on which $\B$ acts as
an automorphism group of an SVOA in \cite{H}.
Our approach is the same as that of \cite{H} and there are
some overlap with this article. 
But our method seems simpler since it is based on explicit
calculations of the intertwining 
operators for the Ising models so that we don't need some assumptions
that are supposed in \cite{H}.

Let $e$ be a rational conformal vector with c.c. $\hf$ in $V^\nat_\R$.
(The existence of such a vector is clear.)
It is shown in \cite{M3} that the Ising model VOA $L(\hf,0)_\R$ over
$\R$ is also rational so that we can apply Theorem \ref{cutting SVOA}
to $V^\nat_\R$.
By Theorem \ref{cutting SVOA}, we can obtain a simple SVOA $T_e^\nat
(0)\oplus T_e^\nat (\hf )$ from $V^\nat_\R$.
Following H\"{o}hn \cite{H}, we set $\VB_\R =T_e^\nat (0)\oplus
T_e^\nat (\hf )$ and call it the Babymonster SVOA.
In $V^\nat_\R$, for every conformal vector $e$, its
Miyamoto-involution $\tau_e$ is not trivial so that
$V^\nat_e(\fr{1}{16})$ is not zero.  
Therefore we also obtain an irreducible $\Z_2$-twisted
$\VB_\R$-module $T_e^\nat(\fr{1}{16})$. We set $(\VBT )_\R :=
T_e^\nat(\fr{1}{16})$. 
Note that  the algebraic structures of $\VB_\R$ and $(\VBT )_\R$ are
independent of the choice of a conformal vector $e$  since every
conformal vector with c.c. $\hf$ in $V^\nat_\R$ is conjugate under the 
Monster $\M$ so that the structures of them are uniquely determined by
that of $V^\nat_\R$. 

\begin{thm}\label{Hohn}
  (1)\ (\cite{H}) The SVOA $\VB_\R$ obtained from $V^\nat_\R$ by
  cutting off the Ising models is a simple SVOA whose full
  automorphism group contains $2\times \B$.
  \vsb\\
  (2)\ The piece $(\VBT )_\R$ obtained from $V^\nat_\R$ is an
  irreducible $\Z_2$-twisted $\VB_\R$-module.
\end{thm}

\pf
The first half assertion of (1) and (2) follow from Theorem
\ref{cutting SVOA}.  
So it remains to show that $2\times \B$ acts on $\VB_\R$ as an
automorphism group of SVOA.
For any conformal vector $e\in V^\nat_\R$, the corresponding
involution $\tau_e$ is unique so that it follows from $g \tau_e g^{-1}
=\tau_{ge}$ for all $g\in \M$ that every element in $C_{\M}(\tau_e )$
fixes $e$. Hence, $C_{\M}(\tau_e )$ acts on $T_e^\nat (0)$.
Furthermore, $C_{\M}(\tau_e )$ leaves the space
$\psi_{-\hf}\vac\tensor T_e^\nat (\hf )\subset V^\nat_\R$ invariant
since each element of $C_{\M}(\tau_e )$ and the action of $e$ commute
on $V^\nat_\R$.   
Therefore we may think $C_{\M}(\tau_e )$ also acts on $T_e^\nat (\hf
)$. 
Since $C_{\M}(\tau_e )=\lfm \tau_e \rfm\cd \B$ and $\tau_e$ acts on
$T_e^\nat (h)$ trivially for $h=0,\hf$, we have an injection from $\B$ 
into the full automorphism group of the SVOA $\VB_\R =T_e^\nat (0)
\oplus T_e^\nat (\hf )$.
\qed

One can expect that the full automorphism group of $\VB_\R$ is the
Baby monster $\B$. 
However, we could not determine it is true or not in this paper. 
But we can prove that $\aut (\VB_\R )$ is finite. To prove this,
we need some results from the quantum Galois theory for SVOAs.

Let $V$ be a simple SVOA over $\C$ and $G$ be a finite subgroup of
$\aut (V)$. 
For $\chi\in\irr (G)$, we set $V^\chi$ to be the sum of all
irreducible components of $V$ on which $G$ acts as $\chi$.
Clearly, $V^\chi$ has a $\Z_2$-grading and we denote such a $\Z_2$-grade
decomposition by $V^\chi=V^{0,\chi}\oplus V^{1,\chi}$.
Let $M_\chi$ be the irreducible representation of $G$ on which $G$
acts as $\chi$. Then we have a decomposition $V^i_\chi =M_\chi\tensor
\hom_G (M_\chi,V^i)$ as a $G$-module.
Setting $V^i_\chi :=\hom_G(M_\chi,V^i)$ for $i=0,1$ and $\chi\in\irr
(G)$, we have the following decomposition of $V$ as $\C [G]\tensor
V^G$-module. 
$$
  V=\bigoplus_{\chi\in\irr (G)} M_\chi\tensor (V^0_\chi \oplus 
  V^1_\chi ).
$$
Then, as an extension of \cite{DM} and \cite{HMT}, we have

\begin{thm} \label{quantum Galois}
  Under the above setting, the followings hold.
  \vsb\\
  (1)\ $V^\chi\ne 0$ for all $\chi \in \irr (G)$.
  \vsb\\
  (2)\ If $V^G$ is a VOA, then $V^i_\chi$, $i=0,1$, $\chi\in\irr (G)$
  except trivial modules are non-isomorphic irreducible $V^G$-modules. 
  \vsb\\
  (3)\ If $V^G$ is an SVOA, then $V^0_\chi \oplus V^1_\chi$ are
  non-isomorphic irreducible $\Z_2$-graded $V^G$-modules.
  In particular, none of $V^i_\chi$, $i=0,1$, $\chi\in \irr (G)$ is
  zero and $\C [G]\tensor V^G$ forms a dual pair over $V$.
\end{thm}

Since the proof is just rewriting of the original quantum Galois
theory for VOAs, we omit it.

\begin{rem}
  It is case-by-case whether $V^G$ is a VOA or SVOA.
  If $G$ contains a canonical involution for SVOA, then $V^G$ must be
  a VOA.
\end{rem}

The following theorem is an analogy of Theorem 9.2 in \cite{M3}.

\begin{thm}
  $\aut (\VB_\R )$ is finite.
\end{thm}

\pf
Since $V^\nat_\R$ is a framed VOA, $\VB_\R$ is also a framed SVOA by
its construction. 
(See \cite{DGH} for the definition of the framed VOAs).
It follows from Miyamoto's construction of the moonshine module
\cite{M3} and our definition of the Babymonster SVOA that $\VB_\R$ has
an invariant bilinear form which is positive definite on the even part 
$(\VB_\R )^0$ and negative definite on the odd part $(\VB_\R )^1$.
Suppose that there exists a subgroup $G$ of $\aut (\VB_\R )$ of
infinite order.
Since $(\VB_\R)_1=0$, we can apply Theorem 9.1 in \cite{M1} and so
there are finitely many conformal vectors with c.c. $\hf$ in $\VB_R$.
So we may assume that $G$ fixes all conformal vectors with c.c. $\hf$
in $\VB_\R$. Let $\w$ be the Virasoro vector of $\VB_\R$.
Then we can find a set of mutually orthogonal conformal vectors
$e^1, e^2,\dots, e^{47}$ with c.c. $\hf$ in $\VB_\R$ such that $\w
=e^1+\cds +e^{47}$. 
Set $P=\lfm \tau_{e^i}, \sigma | i=1,\cdots, 47\rfm$, where
$\tau_{e^i}$ are Miyamoto's involutions and $\sigma$ is a canonical
involution on $\VB_\R$.
By the definition of $\tau_{e^i}$, $P$ is an elementary abelian
$2$-group. Let $\VB_\R =\oplus_{\chi \in \irr (P)} (\VB_\R )^\chi$ be
the decomposition of $\VB_\R$ into the direct sum of eigenspaces of
$P$. 
Since $G$ fixes all $e^i$, $[G,P]=1$ and hence $G$ leaves all $(\VB
_\R )^\chi$ invariant. In particular, $G$ acts on $(\VB_\R )^P$.
Since $P$ contains a canonical involution $\sigma$, $(\VB_\R )^P$ is a 
VOA and is isomorphic to a code VOA $M_D=\oplus_{\alpha\in D}
M_\alpha$ for some even linear code $D$. See \cite{M2} for the
description of the code VOA $M_D$.
Since $G$ fixes all $e^i$, $G$ fixes all elements in $\ising^{\tensor 
47}_\R \subset M_D$ so that $g\in G$ acts on $M_\alpha$ as a scalar
$\lambda_\alpha (g)$.
Since $(\VB_\R )^0$ has a positive definite invariant form, we have
$\lambda_\alpha (g)=\pm 1$ for all $g\in G$.
Therefore, by taking finite index, we may assume that $G$ fixes all
elements in $(\VB_\R )^P$. 
By considering the complexification if necessary, we see that $(\VB_\R
)^\chi$ is an irreducible $(\VB_\R )^P$-module by Theorem \ref{quantum
Galois} so that $g\in G$ acts on $(\VB_\R )^\chi$ as a scalar $\mu_\chi
(g)$. By the same arguments as above, we have a contradiction even
though if we make a complexification. Hence, $\aut (\VB_\R )$ is finite.
\qed

\end{document}